\newcount\secno
\newcount\prmno
\newif\ifnotfound
\newif\iffound

\def\namedef#1{\expandafter\def\csname #1\endcsname}
\def\nameuse#1{\csname #1\endcsname}

\long\def\ifundefined#1#2#3{\expandafter\ifx\csname
  #1\endcsname\relax#2\else#3\fi}
\def\hwrite#1#2{{\let\the=0\edef\next{\write#1{#2}}\next}}

\toksdef\ta=0 \toksdef\tb=2
\long\def\leftappenditem#1\to#2{\ta={\\{#1}}\tb=\expandafter{#2}%
                                \edef#2{\the\ta\the\tb}}
\long\def\rightappenditem#1\to#2{\ta={\\{#1}}\tb=\expandafter{#2}%
                                \edef#2{\the\tb\the\ta}}

\def\lop#1\to#2{\expandafter\lopoff#1\lopoff#1#2}
\long\def\lopoff\\#1#2\lopoff#3#4{\def#4{#1}\def#3{#2}}

\def\ismember#1\of#2{\foundfalse{\let\given=#1%
    \def\\##1{\def\next{##1}%
    \ifx\next\given{\global\foundtrue}\fi}#2}}

\def\section#1{\vskip1truecm
               \global\def\currenvir{section}
               \global\advance\secno by1\global\prmno=0
               {\bf \number\secno. {#1}}
               }

\def\subsection{\global\def\currenvir{subsection}
                \global\advance\prmno by1
                \smallskip \ind{ (\number\secno.\number\prmno) }}
\def\subsec{\global\def\currenvir{subsection}
                \global\advance\prmno by1
                { (\number\secno.\number\prmno)\ }}

\def\proclaim#1{\global\advance\prmno by 1
                {\bf #1 \the\secno.\the\prmno$.-$ }}

\long\def\th#1 \enonce#2\endth{%
   \medbreak\proclaim{#1}{\it #2}\global\def\currenvir{th}\smallskip}

\def\bib#1{\rm #1}
\long\def\thr#1\bib#2\enonce#3\endth{%
\medbreak{\global\advance\prmno by 1\bf#1\the\secno.\the\prmno\ 
\bib{#2}$\!.-$ } {\it
#3}\global\def\currenvir{th}\smallskip}
\def\rem#1{\global\advance\prmno by 1
\par{\it #1} \the\secno.\the\prmno$.-$ }


\def\isinlabellist#1\of#2{\notfoundtrue%
   {\def\given{#1}%
    \def\\##1{\def\next{##1}%
    \lop\next\to\za\lop\next\to\zb%
    \ifx\za\given{\zb\global\notfoundfalse}\fi}#2}%
    \ifnotfound{\immediate\write16%
                 {Warning - [Page \the\pageno] {#1} No reference found}}%
                \fi}%
\def\ref#1{\ifx\labellist\empty{\immediate\write16
                 {Warning - No references found at all.}}
               \else{\isinlabellist{#1}\of\labellist}\fi}

\def\newlabel#1#2{\rightappenditem{\\{#1}\\{#2}}\to\labellist}
\def\labellist{}

\def\label#1{%
  \def\given{th}%
  \ifx\given\currenvir%
    {\hwrite\lbl{\string\newlabel{#1}{\number\secno.\number\prmno}}}\fi%
  \def\given{section}%
  \ifx\given\currenvir%
    {\hwrite\lbl{\string\newlabel{#1}{\number\secno}}}\fi%
  \def\given{subsection}%
  \ifx\given\currenvir%
    {\hwrite\lbl{\string\newlabel{#1}{\number\secno.\number\prmno}}}\fi%
  \def\given{subsubsection}%
  \ifx\given\currenvir%
  {\hwrite\lbl{\string%
    \newlabel{#1}{\number\secno.\number\subsecno.\number\subsubsecno}}}\fi
  \ignorespaces}

\newwrite\lbl

\def\openall{\openout\lbl=\jobname.lbl}

\newread\testfile
\def\lookatfile#1{\openin\testfile=\jobname.#1
    \ifeof\testfile{\immediate\openout\nameuse{#1}\jobname.#1
                    \write\nameuse{#1}{}
                    \immediate\closeout\nameuse{#1}}\fi%
    \immediate\closein\testfile}%

\def\begin{\newlabel{K3}{1.3}
\newlabel{gen}{1.4}
\newlabel{diag}{1.5}
\newlabel{eclat}{1.6}
\newlabel{Fano}{1.7}
\newlabel{bog}{2.2}
\newlabel{cor}{2.3}
\newlabel{flop}{2.6}
\newlabel{\delta }{3.2}
\newlabel{main}{3.3}
\newlabel{sh}{3.4}
\newlabel{j*}{3.7}
\newlabel{X3}{3.8}
\newlabel{dh}{3.9}
\newlabel{io}{3.10}}

\magnification 1250
\pretolerance=500 \tolerance=1000  \brokenpenalty=5000
\mathcode`A="7041 \mathcode`B="7042 \mathcode`C="7043
\mathcode`D="7044 \mathcode`E="7045 \mathcode`F="7046
\mathcode`G="7047 \mathcode`H="7048 \mathcode`I="7049
\mathcode`J="704A \mathcode`K="704B \mathcode`L="704C
\mathcode`M="704D \mathcode`N="704E \mathcode`O="704F
\mathcode`P="7050 \mathcode`Q="7051 \mathcode`R="7052
\mathcode`S="7053 \mathcode`T="7054 \mathcode`U="7055
\mathcode`V="7056 \mathcode`W="7057 \mathcode`X="7058
\mathcode`Y="7059 \mathcode`Z="705A
\def\spacedmath#1{\def\packedmath##1${\bgroup\mathsurround =0pt##1\egroup$}
\mathsurround#1
\everymath={\packedmath}\everydisplay={\mathsurround=0pt}}
\def\nospacedmath{\mathsurround=0pt
\everymath={}\everydisplay={} } \spacedmath{2pt}

\def\phfl#1#2{\normalbaselines{\baselineskip=0pt
\lineskip=10truept\lineskiplimit=1truept}\nospacedmath\smash 
{\mathop{\hbox to 8truemm{\rightarrowfill}}
\limits^{\scriptstyle#1}_{\scriptstyle#2}}}
\def\hfl#1#2{\normalbaselines{\baselineskip=0truept
\lineskip=10truept\lineskiplimit=1truept}\nospacedmath\smash
{\mathop{\hbox to
12truemm{\rightarrowfill}}\limits^{\scriptstyle#1}_{\scriptstyle#2}}}
\def\diagram#1{\def\normalbaselines{\baselineskip=0truept
\lineskip=10truept\lineskiplimit=1truept}   \matrix{#1}}
\def\vfl#1#2{\llap{$\scriptstyle#1$}\left\downarrow\vbox to
6truemm{}\right.\rlap{$\scriptstyle#2$}}
\def\ufl#1#2{\llap{$\scriptstyle#1$}\left\uparrow\vbox
to 6truemm{}\right.\rlap{$\scriptstyle#2$}}
\def\mono{\lhook\joinrel\mathrel{\longrightarrow}}
\def\iso{\vbox{\hbox to .8cm{\hfill{$\scriptstyle\sim$}\hfill}
\nointerlineskip\hbox to .8cm{{\hfill$\longrightarrow $\hfill}} }}

\def\sdir_#1^#2{\mathrel{\mathop{\kern0pt\oplus}
\limits_{#1}^{#2}}}
\def\pprod_#1^#2{\raise
2pt \hbox{$\mathrel{\scriptstyle\mathop{\kern0pt\prod}
\limits_{#1}^{#2}}$}}

\font\eightrm=cmr8         \font\eighti=cmmi8
\font\eightsy=cmsy8        \font\eightbf=cmbx8
\font\eighttt=cmtt8        \font\eightit=cmti8
\font\eightsl=cmsl8        \font\sixrm=cmr6
\font\sixi=cmmi6           \font\sixsy=cmsy6
\font\sixbf=cmbx6\catcode`\@=11
\def\eightpoint{%
  \textfont0=\eightrm \scriptfont0=\sixrm \scriptscriptfont0=\fiverm
  \def\rm{\fam\z@\eightrm}%
  \textfont1=\eighti  \scriptfont1=\sixi  \scriptscriptfont1=\fivei
  \def\oldstyle{\fam\@ne\eighti}\let\old=\oldstyle
  \textfont2=\eightsy \scriptfont2=\sixsy \scriptscriptfont2=\fivesy
  \textfont\itfam=\eightit
  \def\it{\fam\itfam\eightit}%
  \textfont\slfam=\eightsl
  \def\sl{\fam\slfam\eightsl}%
  \textfont\bffam=\eightbf \scriptfont\bffam=\sixbf
  \scriptscriptfont\bffam=\fivebf
  \def\bf{\fam\bffam\eightbf}%
  \textfont\ttfam=\eighttt
  \def\tt{\fam\ttfam\eighttt}%
  \abovedisplayskip=9pt plus 3pt minus 9pt
  \belowdisplayskip=\abovedisplayskip
  \abovedisplayshortskip=0pt plus 3pt
  \belowdisplayshortskip=3pt plus 3pt 
  \smallskipamount=2pt plus 1pt minus 1pt
  \medskipamount=4pt plus 2pt minus 1pt
  \bigskipamount=9pt plus 3pt minus 3pt
  \normalbaselineskip=9pt
  \setbox\strutbox=\hbox{\vrule height7pt depth2pt width0pt}%
  \normalbaselines\rm}\catcode`\@=12

\newcount\noteno
\noteno=0
\def\up#1{\raise 1ex\hbox{\sevenrm#1}}
\def\note#1{\global\advance\noteno by1
\footnote{\parindent0.4cm\up{\number\noteno}\
}{\vtop{\eightpoint\baselineskip12pt\hsize15.5truecm\noindent
#1}}\parindent 0cm}
\font\san=cmssdc10

\def\sym{\hbox{\san \char83}}

\def\pc#1{\tenrm#1\sevenrm}
\def\tx{\kern-1.5pt -}
\def\cqfd{\kern 2truemm\unskip\penalty 500\vrule height 4pt 
depth 0pt width 4pt\medbreak} 
\def\carre{\kern 2truemm\unskip\penalty 500\vrule height  4pt
depth 0pt width 4pt}
\def\virg{\raise
.4ex\hbox{,}}
\def\decale#1{\smallbreak\hskip 28pt\llap{#1}\kern 5pt}
\def\no{n\up{o}\kern 2pt}
\def\ind{\par\hskip 1truecm\relax}
\def\indp{\par\hskip 0.5truecm\relax}

\def\rond{\kern 1pt{\scriptstyle\circ}\kern 1pt}

\def\Ker{\mathop{\rm Ker}\nolimits}

\def\dim{\mathop{\rm dim}\nolimits}

\def\pv{\smallskip {\it Proof} :\ }
\def\pr{\mathop{\rm pr}\nolimits}
\def\br#1{\langle #1\rangle} 
\def\ss{S^{\{2\}_{}}}
\def\sh{S^{[2]}_{}}
\def\bt{\raise-0pt\hbox{$\ \scriptstyle\boxtimes\ $}}
\def\ti{.\thinspace}
\def\du{^{\scriptscriptstyle\vee}}
\def\Q{\Bbb Q}
\def\P{\Bbb P}
\def\C{\Bbb C}

\def\Z{\Bbb Z}

\frenchspacing
\input amssym.def
\input amssym
\vsize = 25truecm
\hsize = 16truecm
\voffset = -.5truecm
\parindent=0cm
\baselineskip15pt
\input xy
\xyoption{all}

\begin
\centerline{\bf On the splitting of the Bloch-Beilinson filtration}
\smallskip
\smallskip \centerline{Arnaud {\pc BEAUVILLE}} 
\vskip1.2cm

{\bf Introduction} \smallskip 
\ind This paper deals with the {\it Chow ring}  $CH(X)$ (with rational
coefficients) of a smooth projective variety $X$   -- that is, 
the $\Q$\tx algebra of algebraic cycles  on
$X$, modulo rational equivalence. This is a basic invariant of
the variety
$X$, which may be thought of as an algebraic counterpart of the
cohomology ring of a compact manifold; in fact there is a $\Q$\tx
algebra homomorphism
$c^{}_X:CH(X)\rightarrow H(X,\Q)$, the {\it cycle class map}. But
unlike the cohomology ring, the Chow ring, and in particular the
kernel of $c^{}_X$, is  poorly understood.
\ind Still some  insight into the structure of this ring is
provided  by the deep conjectures of Bloch and Beilinson. They
predict the existence of a functorial ring filtration
$(F^j)_{j\ge 0}$ of $CH(X)$, with
$CH^p(X)=F^0CH^p(X)\supset\ldots \supset F^{p+1}(X)=0 $ and  
$F^1CH(X)=\Ker c^{}_X$.
We refer to [J] for a  discussion of the various candidates for such a 
filtration and the consequences of its existence.
\ind The existence of that filtration is not even known for
 an abelian variety $A$. In that case, however, there
is a canonical  {\it ring graduation} given by
$CH^p(A)=\sdir_{s}^{}CH^p_s(A)$, where $CH^p_s(A)$ is the
subspace of elements $\alpha \in CH^p(A)$ with $
k_A^*\alpha =k^{2p-s}\alpha $ for all $k\in\Z$ ($k_A$
denotes the endomorphism
$a\mapsto ka$ of $A$) [B2]. 
 Unfortunately this does not define the required filtration
because the vanishing of the terms $CH^p_s(A)$ for $s<0$ is not
known in general -- in fact, this vanishing is essentially equivalent to
the existence of the Bloch-Beilinson filtration. So if the
Bloch-Beilinson filtration indeed exists, it {\it splits} in the sense
that it is the filtration associated to a graduation of $CH(A)$.
\ind In [B-V] we observed that this also  happens  for a K3
surface $S$. Here the filtration is essentially trivial; the fact that
it splits means that the image of the intersection product
$CH^1(S)\otimes CH^1(S)\rightarrow CH^2(S)$ is
always one-dimensional -- an easy but somewhat surprising
property.
\ind The motivation for this paper was to understand whether the
splitting of the Bloch-Beilinson filtration for abelian varieties and
K3 surfaces is accidental or part of a more general framework.
Now asking for a conjectural splitting
of a conjectural filtration may look like a rather idle occupation.
The point we want to make is that the mere existence of such a
splitting has quite concrete consequences, which at least in some
cases can be tested. We will restrict for simplicity to the case of
regular varieties, that is, varieties $X$ for which 
 $F^1CH^1(X)=0$. Then  if the
filtration comes from a graduation, any product of divisors must
have degree $0$; therefore, if we denote by $DCH(X)$ the
sub-algebra  of
$CH(X)$ spanned by divisor classes, {\it the cycle class
map$$c^{}_X:DCH(X)\longrightarrow H(X)$$is injective}. In
other  words, any polynomial relation $P(D_1,\ldots ,D_s)=0$
between divisor classes which hold in cohomology must hold in
$CH(X)$. We will call this property the {\it weak splitting property}.
Despite its name it is rather restrictive: it implies for instance the
existence of a class $\xi_X\in CH^n(X)$, with $n=\dim X$, such that
$$D_1\cdot \ldots \cdot D_n=\deg(D_1\cdot \ldots \cdot D_n)
\cdot \xi_X\quad\hbox{in }CH^n(X)$$for any divisor classes
$D_1,\ldots ,D_n$ in $CH^1(X)$.
\ind What kind of varieties can we expect to have the weak splitting
 property? A natural class containing  abelian varieties and
K3 surfaces is that of Calabi-Yau
varieties, but that turns out to be too optimistic -- it is quite easy to
give counter-examples (Example \ref{Fano}\ti b)). A more restricted
class is that of holomorphic symplectic manifolds -- projective
manifolds admitting an everywhere non-degenerate
holomorphic 2-form. We want to propose the following conjecture:
\smallskip 
{\bf Conjecture}$.-$ {\it A symplectic {\rm (}projective{\rm )}
manifold satisfies the weak splitting property}.
\ind We have to admit that the evidence we are able to provide is
not overwhelming. We will prove that the weak splitting property is
invariant under some simple birational transformations
called Mukai flops (Proposition \ref{flop}). We will also prove that
the conjecture holds for the simplest examples of symplectic
manifolds, the Hilbert schemes $\sh$ and $S^{[3]}_{}$ associated to
a K3 surface $S$ (Proposition
\ref{main}). Already for  $S^{[3]}_{}$ the proof is intricate, and
makes use of some nontrivial relations in the Chow rings of $S^2$
and
$S^3$ established in [B-V]. We hope that this might indicate
 a deep connection between the symplectic structure and the
  Bloch-Beilinson filtration, but we have not even a conjectural
formulation of what this connection could be. 

\section{Intersection of divisors}
\subsection Let $X$ be a  projective (complex) manifold.  We
denote by $CH(X)$ and $H(X)$ the Chow and cohomology rings with
rational coefficients, and by  $CH(X,\C)$ and
$H(X,\C)$ the corresponding rings with complex coefficients.
We  denote by
$DCH(X)$  the sub-algebra of $CH(X)$ spanned by divisor classes.
We will say that $X$ has the {\it weak splitting property}  if the
cycle class map $c^{}_X:DCH(X)\rightarrow CH(X)$ is injective.

 \rem{Remark} The property as stated implies that
$CH^1(X)$ is finite-dimensional, that is,  $X$ is  regular in the
sense that $H^1(X,{\cal O}_X)=0$. For irregular varieties the
definition should be adapted, either by considering cycles modulo
algebraic equivalence, or by picking up an appropriate subspace
of
$CH^1(X)$. We will restrict ourselves to regular varieties in
what follows.
\rem{Examples}\label{K3} a) A regular surface $S$
satisfies the weak splitting property if and only if the image of the
intersection map
$CH^1(S)\otimes CH^1(S)\rightarrow CH^2(S)$ has rank 1;
in other words, there exists a class $\xi_S\in CH^2(S)$, of degree
1, such that $C\cdot D=\deg(C.D)\,\xi_S$ for all curves $C,D$
on $S$. This is the case when $S$ is a K3 surface, or also an 
 elliptic surface over $\P^1$ with a section [B-V].
\ind b) Let $S$ be a K3 surface, $p$ a point of $S$
with$[p]\not=\xi_S$  in $CH^2(S)$. Let
$\varepsilon:\widehat{S}\rightarrow S$ be the blowing-up of
$S$ at $p$. The space $DCH^2(\widehat{S})$ is spanned by
$\varepsilon^*\xi_S$ and $[q]$, where $q$ is any point of
$\widehat{S}$ above $p$. Since the pushforward map
$\varepsilon_*:CH^2(\widehat{S})\rightarrow CH^2(S)$ is an
isomorphism, theses classes are linearly independent in
$CH^2(\widehat{S})$, so the map
$c^2_{\widehat{S}}:DCH^2(\widehat{S})\rightarrow
CH^2(\widehat{S})$ is not injective.
\ind Observe that we get a  family of surfaces
parameterized by $p\in S$, for which the weak splitting
property fails generically, but holds when $p$ lies in
the union of countably many subvarieties of the parameter space.
\smallskip 
\ind c) We will give later (\ref{Fano}) examples of Fano
and Calabi-Yau threefolds which  do not satisfy the
weak splitting property. 
\th Proposition
\enonce Let $X$, $Y$ be two smooth projective {\rm (}regular{\rm
)} varieties.
\indp {\rm a)}  We have
$\displaystyle DCH^p(X\times Y)=\sdir_{r+s=p}^{}
\pr_1^*DCH^r(X)\otimes \pr_2^*DCH^s(Y)$. In particular, 
$X\times Y$ satisfies the weak splitting property if and only if $X$
and $Y$ do.
\indp {\rm b)}  Let $f:X\rightarrow Y$ be a surjective map. If
$c^p_X:DCH^p(X)\rightarrow H^{2p}(X)$ is injective, then so is
$c^p_Y:DCH^p(Y)\rightarrow H^{2p}(Y)$.
\endth\label{gen}
\pv a)  We have $CH^1(X\times
Y)=\pr_1^*CH^1(X)\oplus \pr_2^*CH^1(Y)$ since $X$ and $Y$ are
regular; the assertion a) follows at once.
\ind b) follows from  the commutative diagram\vskip-12pt
$$\diagram{DCH^p(X) & \hfl{c^p_X}{} & H^{2p}(X)\cr
\ufl{f^*}{}&&\ufl{f^*}{}\cr
DCH^p(Y) & \hfl{c^p_Y}{} & H^{2p}(Y)}$$\vskip-8pt and the
injectivity of
$f^*:CH^p(Y)\rightarrow CH^p(X)$ (if $h$ is an ample class in
$CH^1(X)$ and $d=\dim X-\dim Y$, we have
$f_*(h^d)=r\cdot 1_Y$, with $r\in\Q^*$,  and
$f_*(h^d\cdot f^*\xi)=$ $r\,\xi\
$ for
$\xi$ in $CH(Y)$).\cqfd

\subsection\label{diag} We now consider the behaviour of the weak
splitting property when the variety $X$ is blown up  
 along a smooth subvariety $B$. We will use the notation
summarized in the following diagram:
$$\diagram{E &\lhook\joinrel\mathrel{\hfl{i}{}}&
\widehat{X}&\cr
\vfl{\eta}{}&&\vfl{}{\varepsilon}&\cr
B &\lhook\joinrel\mathrel{\hfl{j}{}}&
X&\kern-10pt.\cr}\leqno{(\ref{diag})}$$ We denote by
$c$ the codimension of $B$ in $X$ and by $N$ its normal bundle.
\th Lemma
\enonce Let $p$ be an integer. Assume {\rm :}
\indp {\rm (i)}   The cycle class map $c^q_B:DCH^q(B)\rightarrow
H^{2q}(B)$ is injective for $p-c<q<p\,;$
\indp {\rm (ii)}  The Chern
classes $c_i(N)$ belong to
$DCH(B)\,;$
\indp {\rm (iii)}  The map $c^p_X:CH^p(X)\rightarrow
H^{2p}(X)$ restricted to
$DCH^p(X)+j_*DCH^{p-c}(B)$ is injective.
\ind Then the cycle class map
$c^p_{\widehat{X}}:DCH^p(\widehat{X})\rightarrow
H^{2p}(\widehat{X})$ is injective.
\endth\pv\label{eclat}
 The projection $p:E\rightarrow B$ identifies $E$ to
$\P_B(N\du)$. Let
$h\in CH^1(E)$ be the class of the tautological bundle ${\cal
O}_E(1)$; we have $i^*[E]=-h$, and therefore, for $\xi \in CH(X)$,
$[E]^p\cdot \varepsilon ^*\xi =i_*(i^*[E]^{p-1}\cdot i^*\varepsilon
^*\xi )=(-1)^{p-1}i_*(h^{p-1}\cdot \eta ^*j^*\xi )$.
\ind Since $CH^1(\widehat{X})=\varepsilon^*CH^1(X)\oplus
\Q [E]$,  we get
$$\nospacedmath\displaylines{DCH^p(\widehat{X})\ =\
\varepsilon^*DCH^p(X)+ [E]\cdot
\varepsilon^*DCH^{p-1}(X)+\ldots +\Q [E]^p\hfill\cr 
\hfill\i\ \varepsilon^*DCH^p(X)+i_*\eta ^*DCH^{p-1}(B)+i_*(h\cdot
\eta ^*DCH^{p-2}(B))+\ldots +\Q i_*h^{p-1}\ .}$$ 
For $q\ge c$ we have a relation $h^q=h^{c-1}\cdot
\eta ^*c_{q,c-1}+\ldots +\eta ^*c_ {q,0}$, where the $c_{i,j}$ are
polynomial in the Chern classes of $N$;  by our
hypothesis (ii) these classes lie in $DCH(B)$. Moreover the ``key
formula" [F, 6.7]
$$i_*(\gamma\cdot \eta ^*\xi)=\varepsilon^*j_*\xi\qquad
\hbox{for }\xi\in CH(B)\ ,$$with
$\gamma=h^{c-1}+h^{c-2}\cdot \eta ^*c_1(N)+\ldots
+\eta ^*c_{c-1}(N)$, implies $$i_*(h^{c-1}\cdot \eta ^*DCH^{p-c}(B))\i
\varepsilon^*j_*DCH^{p-c}(B)+\sum_{k=0}^{c-2} i_*(h^{k}\cdot
\eta ^*DCH^{p-k-1}(B))\ ,$$so that we finally get
$$DCH^p(\widehat{X})\ \i\ \varepsilon^*\bigl(DCH^p(X)+
j_*DCH^{p-c}(B)\bigr)+\sum_{k=0}^{c-2} i_*(h^{k}\cdot
\eta ^*DCH^{p-k-1}(B))$$
Since the map $$\eqalign{ H^{2p}(X)\oplus\sum_{k=0}^{c-2}
H^{2(p-k-1)}(B)\ & \longrightarrow\  H^{2p}(\widehat{X})\cr
(\alpha\ ;\ \beta_0,\ldots ,\beta_{c-2})\qquad& \longmapsto\ 
\varepsilon^*\alpha+\sum_ki_*(h^k\cdot \eta^*\beta_{k})}$$
is an isomorphism  (see for instance [Jo]), our
hypotheses (i) and (iii) ensure that $c^p_{\widehat{X}}$ is
injective.\cqfd
\rem{Examples}\label{Fano} a) Take $X=\P^3$,
and let $B$ be a smooth curve, of degree $d$ and genus $g$. Let
$\ell $ be the class of a hyperplane in
$\P^3$,
$\ell _B$ its pull back to $B$. The space
$DCH^2(\widehat{X})$ is generated by 
$$\varepsilon^*\ell ^2\quad,\quad \varepsilon^*\ell\cdot
[E]=i_*p^*\ell_{B}
\quad,\quad [E]^2=-i_*h=i_*p^*c_1(N)-\varepsilon^*[B]
$$
We have $c_1(N)=4\ell _B+K_B$, so $DCH^2(\widehat{X})$
contains the elements $i_*p^*\ell_{B}$ and $i_*p^*K_B$. 
\ind The map $i_*p^*:CH^1(B)\rightarrow CH^2(X)$ induces an
isomorphism of  the subspace of degree 0 divisor classes on $B$
onto the subspace of homologically trivial classes in $CH^2(X)$.
If we choose $\ell _B$ non proportional to $K_B$ in $CH^1(B)$,
the class $i_*p^*(d\,K_B-(2g-2)\ell _B)$ in $DCH^2(\widehat{X})$
is homologically trivial, but non-trivial.
Thus the map $c^2_{\widehat{X}}:DCH^2(\widehat{X})\rightarrow
H^4(\widehat{X})$ is not injective.
\ind If  $B$ is a scheme-theoretical intersection of
cubics, 
$\widehat{X}$ is a Fano variety [M-M] -- we can take for instance
$B$ of genus 2 and $\ell _B$ a general divisor class of degree 5
(or $B$ of genus 3 and $\ell _B$  general of degree 6, or $B$ of
genus 5 and $\ell _B\equiv K_B-p$ for $p$ a  general point of $B$).
Note that by making the linear system vary we get again 
families where the general member does not satisfy the weak
splitting property, while countably many special members of the
family do satisfy it.
\smallskip  
\ind b) Going on with the Fano case, let $D$ be a smooth
divisor in $|-2K_X|$, and let $V\rightarrow X$ be the double
covering of $X$ ramified along $D$. Then by the
above example and Proposition
\ref{gen}\ti b),   
$V$ {\it is a Calabi-Yau threefold which does not satisfy the weak
splitting property}. 
 
\section{The weak splitting property for symplectic manifolds}
\subsection By a symplectic manifold  we
mean here a simply-connected projective manifold which admits a 
holomorphic, everywhere
non-degenerate  2-form. The manifold is said to be
{\it irreducible} if the 2-form is unique up to a scalar; any
symplectic  manifold admits a canonical decomposition as a
product of  irreducible ones. In view of Proposition
\ref{gen}\ti a), we may restrict ourselves to irreducible
symplectic manifolds. 
\ind Let $X$ be an irreducible symplectic manifold, of dimension
$2r$.  Recall that the space $H^2(X)$ admits
a canonical quadratic form
$q$ ([B1], [H]) with the following  properties:
\indp -- every class $x\in H^2(X,\C)$ with $q(x)=0$ satisfies
$x^{r+1}=0$;
\indp -- there exists  $\lambda \in \Q$ such that
$\int_Xx^{2r}=\lambda \,q(x)^r$ for all $x\in H^2(X,\C)$, where
$\int_X$ is the canonical isomorphism $H^{2r}(X,\C)\iso \C$. 
\ind In
fact the following more precise statement has been proved by
Bogomolov:
\th Proposition 
\enonce  Let $V$ be a subspace of $H^2(X,\C)$ such
that the restriction of $q$ to $V$ is non-degenerate {\rm (}for
instance
$V=H^2(X,\C)$ or
$V=CH^1(X,\C))$.  The kernel of the map
$\sym V\rightarrow H(X,\C)$ is the ideal of $\sym V$ spanned by 
the elements $x^{r+1}$ for
$x\in V,\ q(x)=0$.
\endth\label{bog}
\pv The case $V=H^2(X,\C)$ is the main result of [Bo], but the proof
given there implies the slightly more general statement \ref{bog}.
Namely, define 
 $A(V)$ as the quotient of $\sym V$ by the ideal 
spanned by the elements $x^{r+1}$ for $x\in V,\ q(x)=0$. Then
Lemma 2.5 in [Bo] says that
 $A(V)$ is a  finite-dimensional graded Gorenstein $\C$\tx algebra,
with socle in degree $2r$ -- in other words,
$A_{2r}(V)$ is  one-dimen\-si\-onal, and  the multiplication
pairing $A_d(V)\times A_{2r-d}(V)\rightarrow
A_{2r}(V)\cong\C$ is a perfect duality. 
\ind Since  any element $x$ of $H^2(X,\C)$ with $q(x)=0$ satisfies
$x^{r+1}=0$, we get a $\C$\tx algebra 
homomorphism $u:A(V)\rightarrow H(X,\C)$. The kernel of $u$ is 
an ideal of $A(V)$; if it is non-zero, it contains the minimal ideal
$A_{2r}(V)$ of $A(V)$. But this is impossible because $V$
contains an element $h$ with $q(h)\not=0$, hence with
$h^{2r}\not=0$.\cqfd
\th Corollary
\enonce The following conditions are equivalent:
\indp {\rm (i)} The cycle class map $c^{}_X:DCH(X)\rightarrow
H(X)$ is injective {\rm (}that is, $X$ satisfies the weak splitting
property{\rm );}
\indp{\rm (ii)} The map $c^{r+1}_X:DCH^{r+1}(X)\rightarrow
H^{2r+2}(X)$ is injective;
\indp{\rm (iii)} Every element $x$ of $CH^1(X,\C)$ with $q(x)=0$
satisfies $x^{r+1}=0$ {\rm (}in $CH^{r+1}(X,\C)).$
\endth\label{cor}\pv Consider the diagram
$$\xymatrix{\sym CH^1(X,\C)\ar[d]_{u }\ar[dr]^{v }\\
DCH(X,\C)\ar[r]_{c^{}_X} & H(X,\C)\ .
}$$
The injectivity of $c$ is equivalent to $\Ker v \i\Ker u  $.
In view of the Proposition, this is exactly condition (iii), and it is
equivalent to $\Ker v^{r+1} \i\Ker u^{r+1}$.\cqfd
\rem{Remark} Assume that there is an element $\alpha\in
CH^1(X)$ with $q(\alpha)=0$ -- this is the case for instance if
$\dim_{\Q} CH^1(X)\ge 5$. Then the set of such elements  is Zariski
dense in the quadric $q=0$ of
$CH^1(X,\C)$. Thus the conditions of the Corollary are also
equivalent to:
\indp ${\rm (iii')}$ {\it Every element $x$ of $CH^1(X)$ with
$q(x)=0$ satisfies} $x^{r+1}=0$.
\ind A possible proof of $({\rm iii'})$ could be as follows. It seems
plausible that  the subset of {\it nef} classes $x\in CH^1(X)$ with
$q(x)=0$ is Zariski dense in
the quadric $q=0$ (this holds at least when $X$ is a K3
surface). If this is  the case, it would be enough to prove
$({\rm iii'})$ for nef classes. Now it is a standard
conjecture  (see  [S]) that  a
nef class $x\in CH^1(X)$ with
$q(x)=0$ should be the pull back of the class of a hyperplane in
$\P^r$ under a Lagrangian fibration
$f:X\rightarrow\P^r$, so that
$x^{r+1}=f^*(h^{r+1})=0$.
\smallskip 
\subsection We will now consider the behaviour of 
the weak splitting property under a Mukai flop. Let $X$ be an
irreducible symplectic manifold, of dimension $2r$; assume that $X$
contains a subvariety
$P$ isomorphic to $\P^r$. Then $P$ is a Lagrangian subvariety,
and its normal bundle in $X$ is isomorphic to
$\Omega^1_P$. We blow up $P$ in $X$, getting our standard
diagram
$$\diagram{E &\lhook\joinrel\mathrel{\hfl{i}{}}&
\widehat{X}&\cr
\vfl{\eta}{}&&\vfl{}{\varepsilon}&\cr
P &\lhook\joinrel\mathrel{\hfl{j}{}}& X&\kern-12pt .\cr}$$

The exceptional divisor $E$ is the cotangent bundle $\P(T_P)$,
which can be identified with the incidence divisor in $
P\times P\du$, where $P\du$ is the 
projective space dual to
$P$. The projection $\eta\du:E\rightarrow P\du$ identifies $E$
to
$\P(T_{P\du})$, and 
$E$ can be blown down to $P\du$ by a map $\varphi
:\widehat{X}\rightarrow X'$, where $X'$ is a 
smooth  algebraic space.  
To remain in our previous framework we will assume
that $X'$ is projective, so that $X'$ is again an irreducible symplectic
manifold. The diagram
$$\xymatrix{&\widehat{X}\ar[dl]_{\varepsilon }\ar[dr]^{\varphi
}&\\X&&X' }$$is called a {\it Mukai flop}. There are many concrete
examples of such  flops, see [M].
\th Proposition
\enonce If $X$ satisfies the weak splitting property, so does $X'$.
\endth\label{flop}\pv Consider the  $\Q$\tx linear map
$\varphi _*\varepsilon^*:CH^1(X)\rightarrow CH^1(X')$. It is
bijective and preserves the canonical
quadratic forms (see e.g. [H], Lemma 2.6). In view of Corollary
\ref{cor}, the Proposition will follow from
\th Lemma
\enonce Let $\alpha\in CH^1(X)$, and $\alpha ':=\varphi _*\varepsilon
^*\alpha $. Then $\alpha '^{r+1}=\varphi _*\varepsilon
^*(\alpha ^{r+1})$.
\endth\pv
We have $\varphi ^*\alpha '=\varepsilon ^*\alpha
+m[E]$ for some $m\in\Q$.   
Let $\ell\in CH^{2r-1} (\widehat{X})$
be the class of  a line contained in a fibre of $\eta\du$; we
have
$\deg([E]\cdot \ell) =-1$, and
$\varepsilon  _*\ell $ is the class of  a line in $P$. Intersecting the
above equality with $\ell $ gives $m=\deg(\alpha _{|P})$, or
equi\-valently $\alpha _{|P}=mk$ in $CH^1(P)$, where $k$ is the
class of a hyperplane in $P$. Then 
$$\varphi ^*\alpha '^{r+1}=(\varepsilon ^*\alpha
+m[E])^{r+1}=\sum_{p=0}^{r+1} {r+1\choose
p}m^{r+1-p}\,\varepsilon ^*\alpha ^p\cdot [E]^{r+1-p}\ .$$
 As in (\ref{eclat}), let $h\in CH^1(E)$ be the class of ${\cal
O}_E(1)$.  For $p\le r$ we have 
 $$\varepsilon ^*\alpha ^p\cdot
[E]^{r+1-p}=(-1)^{r-p}i_*(h^{r-p}\cdot i^*\varepsilon ^*\alpha
^p)=(-1)^{r-p}i_*(h^{r-p}\cdot \eta ^*\alpha_{|P}
^p)\ ,$$\vskip-25pt
$$\varphi ^*\alpha '^{r+1}=\varepsilon
^*\alpha ^{r+1}+m^{r+1}i_*\Bigl(\sum_{p=0}^{r} {r+1\choose
p}(-1)^{r-p}h^{r-p}\,\eta ^*k^p\Bigr)\ .\leqno{\rm Thus}$$
Now since the total Chern class of $T_P$ is $(1+k)^{r+1}$ we
have in $CH^r(E)$
$$\sum_{p=0}^r
{r+1\choose p}(-1)^p h^{r-p}\,\eta ^*k^p=\sum_{p=0}^r
(-1)^ph^{r-p}\,\eta ^*c_p(T_P)=0\ ,$$ 
hence $\varphi ^*\alpha '^{r+1}=\varepsilon ^*\alpha ^{r+1}$.
Applying $\varphi _*$ gives the lemma, hence the
Proposition.\cqfd 
\th Corollary
\enonce Let $X,X'$ be  birationally equivalent projective symplectic
fourfolds. Then $X$ satisfies the weak splitting property if and only
if $X'$ does.
\endth\ind  Indeed any birational map between
projective symplectic fourfolds is a composition of Mukai flops
[W].
\section{The weak splitting property for ${\bf \sh}$ and
${\bf S^{[3]}}$.}
\subsection The simplest symplectic manifolds are K3 surfaces, for
which we have already seen that the weak splitting property holds
(Example \ref{K3}). More precisely [B-V], let $S$ be a K3 surface and
$o$ a point of $S$ lying on a (singular) rational curve $R$. The class
of $o$ in
$CH^2(S)$ is independent of the choice of $R$, and we have, for
every $\alpha ,\beta \in CH^1(S)$,
$$\alpha \cdot \beta =\deg(\alpha \cdot \beta)\,[o]\quad
\hbox{in }\ CH^2(S)\ .$$
\ind Let $\Delta :S\mono S\times S$ be the diagonal embedding. For 
$\alpha \in CH^1(S)$, we have in $CH^3(S\times S)$ ([B-V], Prop.
1.6)\global\advance\prmno by1\label{\delta }
$$\Delta _*\alpha =\pr_1^*\alpha\cdot
\pr_2^*[o]+\pr_1^*[o]\cdot \pr_2^*\alpha \ .\leqno{
(\number\secno.\number\prmno)}$$

 \ind K3 surfaces are  the first instance of a famous series of
symplectic manifolds, the {\it Hilbert schemes} $S^{[r]}$
parameterizing finite subschemes of length $r$ on the K3 surface
$S$.

\th Proposition
\enonce Let $S$ be a $K3$ surface. The symplectic
varieties $\sh$ and $S^{[3]}$ satisfy the weak splitting property.
\endth\label{main}
\pv\kern-10pt\subsec  Let us warm up with the easy case of
$\sh$. Let $\ss$ be the variety obtained by blowing up the
diagonal  of
$S\times S$. The Hilbert scheme
$\sh$ is the quotient of $\ss$ by the
involution which exchanges the factors. In view of Corollary
\ref{bog} and Proposition
\ref{gen}\ti b) it suffices to prove that the cycle class map
$c^3_{\ss}:DCH^3(\ss)\rightarrow H^6(\ss)$ is injective. We will
check that the hypotheses of Lemma \ref{eclat} are satisfied.
Condition (i) is the weak splitting property for $S$. The normal
bundle to the diagonal in $S\times S$
 is $T_S$, so (ii)  means that the class $c_2(T_S)$ belongs to
$DCH^2(S)$; this is proved in ([B-V], thm. 1 c).  Formula (\ref{\delta
}) implies
$\Delta_*CH^1(S)\i DCH^3(S\times S)$, so condition (iii) reduces to
the  injectivity of $c^3_{S\times S}$, which follows from 
Proposition \ref{gen}\ti a) and the corresponding result for
$S$.\cqfd\label{sh}
\subsection Let us pass to the more difficult case of $S^{[3]}$. 
 The Hilbert scheme $S^{[3]}$ is dominated by the {\it
nested Hilbert scheme} $S^{[2,3]}$ which parameterizes pairs
$(Z,Z')\in S^{[2]} \times S^{[3]}$ with $Z\i Z'$; it is isomorphic
to the blow-up of $S\times S^{[2]}$ along the incidence
subvariety ${\cal I}=\{(x,Z)\ |\ x\in Z\}$. Let $\pi :\ss\rightarrow
\sh$ be the quotient map, and $p:\ss\rightarrow S$ the first
projection. Then the map  $j=(p,\pi ):\ss\mono
S\times S^{[2]}$ induces an isomorphism of $\ss$ onto ${\cal
I}$ (see for instance [L], 1.2).  

\ind To prove the theorem, it suffices, by Corollary \ref{cor} and
Proposition \ref{gen}\ti b), to prove that the cycle class map
$DCH^4(S^{[2,3]})\rightarrow H^8(S^{[2,3]})$ is injective. We
will again check that the hypotheses of Lemma
\ref{eclat}  are satisfied. Condition (i) is the injectivity of the cycle
class map $c^3_{\ss}:DCH^3(\ss)\rightarrow H^6(\ss)$, which has
just been proved. Let
$N$ be the normal bundle to the embedding $j:\ss\mono
S\times S^{[2]}$, and $E\i\ss$ the exceptional divisor, which is
the ramification locus of $\pi $. From the exact sequences
$$\nospacedmath \displaylines{0\rightarrow N\du\longrightarrow
p^*\Omega^1_S\oplus
\pi ^*\Omega^1_{\sh}\longrightarrow \Omega^1_{\ss}\rightarrow
0\cr
0\rightarrow \pi ^*\Omega^1_{\sh}\longrightarrow
\Omega^1_{\ss}\longrightarrow {\cal O}_E(-E)\rightarrow
0\cr
0\rightarrow {\cal O}_{\ss}(-2E)\longrightarrow {\cal
O}_{\ss}(-E)\longrightarrow {\cal O}_E(-E)\rightarrow 0}$$we
obtain the equality in K-theory
$[N\du]=[p^*\Omega^1_S]+[{\cal O}_{\ss}(-2E)]-[{\cal
O}_{\ss}(-E)]$. We conclude that $c_2(N)=c_2(N\du)$ belongs to
$DCH^2(\ss)$, so that condition (ii) holds.
\subsection The rest of the proof will be devoted to check condition
(iii), namely the injectivity of
$$DCH^4(S\times \sh)+j_*DCH^2(\ss)\longrightarrow H^8(S\times
\sh)\ .$$
\ind Let us fix some notation. We will use our standard diagram
(\ref{diag})
$$\diagram{E &\lhook\joinrel\mathrel{\hfl{i}{}}&
S^{\{2\}}&\cr
\vfl{\eta}{}&&\vfl{}{\varepsilon}\cr
S &\lhook\joinrel\mathrel{\hfl{\Delta}{}}& S^2&\kern-18pt .\cr}$$
We denote by $p$ and $q$  the two projections of $\ss$ onto $S$.
\ind We define an injective $\Q$\tx linear map
$\iota:CH(S)\rightarrow CH(\sh)$ by $\iota (\xi ):=\pi _*p^*\xi $;
we will use the same notation for cohomology classes. We have
$\pi ^*\iota (\xi )=p^*\xi +q^*\xi $ for $\xi$ in $CH(S)$ or
$H(S)$. Finally if $\alpha\in CH(S)$ and $\xi\in CH(\sh)$ we put 
$\alpha\bt\xi:=\pr_1^*\alpha\otimes \pr_2^*\xi$. 
\ind We have $CH^1(\ss)=p^*CH^1(S)\oplus q^*CH^1(S)\oplus
\Q[E]$.
  In $CH^2(\ss)$ we have  $[E]\cdot
\varepsilon^*\alpha=i_*\eta^*\Delta^*\alpha$ for $\alpha\in
CH^1(S^2)$, and $[E]^2=-\varepsilon^*[\Delta(S)]$. Therefore:
$$DCH^2(\ss)=\Q\,  p^*[o]+\Q\,  q^*[o]\
+\ p^*CH^1(S)\otimes q^*CH^1(S)\ +\
i_*\eta^*CH^1(S)+\Q\,  \varepsilon^*[\Delta(S)]\ .$$
We want to describe the space $j_*DCH^2(S)+DCH^4(S\times \ss)$.

\th Lemma
\enonce Let $\alpha,\beta \in CH^1(S)$. The classes
$j_*p^*[o]$, $j_*(\ p^*\alpha\cdot  q^*\beta )$,  and $j_*i_*\eta
^*\alpha $ belong to
$DCH^4(S\times
\ss)+\Q\, ([o]\bt \iota
([o]))$.\label{j*}
\endth\pv  Let $j':\ss \mono S\times \ss$ be the embedding given
by
$j'(z)=(p(z),z)$, so that $j=(1,\pi )\rond j'$. From the cartesian
diagram
$$\diagram{\ss& \lhook\joinrel\mathrel{\hfl{j'}{}} & S\times
\ss\cr
\vfl{p}{}&&\kern-9pt \vfl{}{(1,p)}\cr
S& \lhook\joinrel\mathrel{\hfl{\Delta}{}} & \kern-9pt S\times
S}\leqno(\ref{j*})$$ we obtain
$j'_*p^*[o]=(1,p)^*\Delta_*[o]=[o]\bt p^*[o]$, hence
$j_*p^*[o]=$ $[o]\bt \iota ([o])$. In the same way we have
$j'_*p^*\alpha=(1,p)^*\Delta_*\alpha$, hence, using (\ref{\delta
}), 
$$j'_*p^*\alpha=\alpha\bt p^*[o]+[o]\bt p^*\alpha\
.$$
Multiplying by
$\pr_2^*q^*\beta$ and using
$\pr_2\rond j'={\rm Id}$ we obtain
$$j'_*(p^*\alpha\cdot q^*\beta)=\alpha\bt (p^*[o]\cdot
q^*\beta)+[o]\bt (p^*\alpha\cdot q^*\beta)\ ,$$
$$j_*(p^*\alpha\cdot q^*\beta)=\alpha\bt \pi _*(p^*[o]\cdot
q^*\beta)+[o]\bt \pi _*(p^*\alpha\cdot q^*\beta)\
.\leqno{\hbox{hence}}$$ 
 For $\alpha ,\beta \in CH^1(S)$,
 put $\br{\alpha ,\beta }:=\deg(\alpha \cdot \beta )$. Then
$$\eqalign{\pi ^*\pi
_*(p^*\alpha\cdot q^*\beta)&=p^*\alpha\cdot
q^*\beta+p^*\beta\cdot q^*\alpha\cr
&=(p^*\alpha+q^*\alpha)(p^*\beta+q^*\beta)-
\br{\alpha,\beta}(p^*[o]+q^*[o])\cr
&=\pi^*\bigl(\iota(\alpha)\iota(\beta)-
\br{\alpha,\beta}\iota([o])\bigr)\ ;}$$we find similarly  $
\pi ^*\pi _*(p^*[o]\cdot q^*\beta)=\pi ^*\iota([o])\iota(\beta)$,
and finally
$$j_*(p^*\alpha \cdot q^*\beta )=\alpha \bt \iota ([o])\iota
([\beta ])\ +\ [o]\bt \bigl(\iota (\alpha )\iota (\beta )-\br{\alpha
,\beta }\,\iota ([o])\bigr)\ .$$
\ind  Let $\gamma \in CH^1(S)$. We have
$\ \iota(\beta )^2\cdot
\iota(\gamma )=\br{\beta ^2}\,\iota([o])\iota(\gamma )+2\br {\beta
\cdot\gamma }\,\iota([o])\iota(\beta )\ $ (this is easily checked
 by applying $\pi^* $  as above). If
$\br{\beta ^2}\not= 0$ we conclude by taking $\gamma  =\beta 
$ that
$\iota([o])\iota(\beta )$ is proportional to
$\iota(\beta )^3$. If
$\br{\beta ^2}= 0$ we can choose $\gamma $ so that $(\beta \cdot
\gamma )\not= 0$; then  $\iota([o])\iota(\beta )$ is proportional
to
$\iota(\beta )^2\iota(\gamma )$. In each case we see that
 $\iota([o])\iota(\beta )$ belongs to $DCH^3(\sh)$, hence the
assertion of the lemma about
$j_*(\ p^*\alpha\cdot  q^*\beta )$. 
\ind Consider finally the cartesian diagram
$$\diagram{E& \lhook\joinrel\mathrel{\hfl{k}{}} & S\times
E\cr
\vfl{\eta}{}&& \vfl{}{(1,\eta)}\cr
S& \lhook\joinrel\mathrel{\hfl{\Delta}{}} & S\times S}$$
with $k(e)=(\eta(e),e)$. Using again (\ref{\delta })
we get $$k_*\eta^*\alpha= (1,\eta)^*\Delta_*\alpha=\alpha\bt 
\eta^*[o]+[o]\bt \eta^*\alpha\ .$$ Pushing
forward in $S\times \sh$ we obtain
$j_*i_*\eta^*\alpha=\alpha\bt  i'_*\eta ^*[o]+[o]\bt
i'_*\eta ^*\alpha$, where $ i'=\pi \rond i$ is the
embedding of $E$ in $\sh$.
\ind  To avoid
confusion let us denote by $\bar E$ the image of $E$ in $\sh$, so
that $\pi ^*[\bar E]=2[E]$. We have $ i'_*\eta ^*\alpha=\pi
_*([E]\cdot p^*\alpha) ={1\over 2}[\bar E]\cdot \iota (\alpha )\in
DCH^2(\sh)$. Likewise $[E]^3=i_*h^2=-24 i_*\eta^*[o]$, hence
$ i'_*\eta ^*[o]=-{1\over 96}\,[\bar E]^3\in
DCH^2(\sh)$. This finishes the proof of the lemma.\cqfd

\ind The lemma and the formula for $DCH^2(\ss)$  show that $
j_*DCH^2(\ss)$ is spanned  modulo $DCH^4(S\times \sh)$ by the
classes
$$ [o]\bt \iota([o])\ ,\quad j_*q^*[o]\
,\quad j_*\varepsilon^*[\Delta(S)]\ .$$
In fact there is one more relation, much more subtle, between
these classes modulo $DCH^4(S\times \sh)$. 

\th Lemma
\enonce We have $$2[o]\bt \iota ([o])-2j_*q^*[o]+ j_*\varepsilon
^*[\Delta(S)]\ \in\ DCH^4(S\times
\sh)\ .$$
\endth\label{X3}\pv
We start from a relation in $CH^4(S^3)$, proved in [B-V, Prop.
3.2]. For $1\le i<j\le 3$, let us denote by $p_{ij}:S^3\rightarrow
S^2$ the projection onto the $i$\tx th and $j$\tx th factors, and by
$p_i:S^3\rightarrow S$ the projection onto the $i$\tx th factor. 
We will write simply $\Delta$  for the diagonal $\Delta(S)\i S^2$,
and $\delta \i S^3$ for the small diagonal, that is, the subvariety of
triples
$(x,x,x)$ for
$x\in S$. Then:
$$[\delta ]- \sum_{i<j, k\not=i,j}p_{ij}^*[\Delta ]\cdot
p_k^*[o]+\sum_{i<j}p_i^*[o]\cdot p_j^*[o]=0\ .$$ Pull back
this relation by the map $\varepsilon^{}_S=(1_S,\varepsilon ):
S\times
\ss\rightarrow S\times S^2$. Since
$$\nospacedmath\displaylines{p_1\rond
\varepsilon^{} _S=\pr_1\ ,\quad  p_2\rond
\varepsilon^{} _S=p\rond \pr_2\ ,\quad p_3\rond \varepsilon^{}
_S=q\rond
\pr_2\ ,\quad p_{23}\rond \varepsilon^{} _S=\varepsilon \ ,\cr
\hbox{we obtain }\hfill\quad \varepsilon _S^*[\delta
]=j'_*\varepsilon ^*[\Delta ]\    ,\qquad 
\varepsilon  _S^*(p_1^*[o]\cdot p_{23}^*[\Delta
])=[o]\bt\varepsilon ^*[\Delta ]\ , \hfill\cr
\varepsilon_S^*(p_2^*[o]\cdot p_3^*[o])=1\bt p^*[o]\cdot
q^*[o] \  ,\quad 
\varepsilon_S^*(p_1^*[o]\cdot p_2^*[o])=[o]\bt p^*[o]\ ,\cr   
\varepsilon_S^*(p_1^*[o]\cdot p_3^*[o])=[o]\bt q^*[o]\ .
 }$$ 
\ind We have $p_{12}\rond \varepsilon _S=(1,p)$, hence  
$\ \varepsilon_S^* p_{12}^*[\Delta ]=j'_*1 $ (see diagram \ref{j*})
and
$ \varepsilon  _S^*(p_3^*[o]\cdot p_{12}^*[\Delta])=j'_*q^*[o]$.
Let $j''=(q,1):\ss\rightarrow S\times \ss$; the same
argument gives
$\ \varepsilon_S^* p_{13}^*[\Delta ]=j''_*1\ $ and
$\ \varepsilon  _S^*(p_2^*[o]\cdot p_{13}^*[\Delta])=j''_*p^*[o]\
$. Finally we have \break$j'_*q^*[o]+j''_*p^*[o]=\pi _S^*j_*q^*[o]$.
Pushing forward by $\pi _S$ we obtain in $CH^4(S\times \sh)$:
$$j_*\varepsilon ^*[\Delta ]-2j_*q^*[o]-[o]\bt \pi _*\varepsilon
^*[\Delta ]+2[o]\bt \iota ([o])+1\bt \iota ([o])^2=0$$
  It remains to observe that  $[o]\bt \pi _*\varepsilon
^*[\Delta ]$ and $1\bt \iota ([o])^2$ belong to
$DCH^4(S\times
\sh)$. Indeed from $[E]^2=-\varepsilon ^*[\Delta ]$ we  deduce
$\pi
_*\varepsilon^*[\Delta]=-{1\over 2}[\bar E]^2\in  DCH^2(\sh)$.
And if $h$ is any element of $CH^1(S)$ with $h^2=d\not=0$, we
have $$\pi ^*\iota (h)^4=6d^2p^*[o]\cdot q^*[o]=3d^2\pi ^*\iota
([o])^2\ ,\ \hbox{ hence }\ \iota ([o])^2\in DCH^4(\sh)\ .\carre$$
\medskip
\subsection\label{dh} For a smooth projective variety $X$, let us
denote by
$DH(X)$ the (graded) subspace of $H(X)$ spanned by intersection
of divisor classes -- that is, the image of $DCH(X)$ by the cycle class
map. It remains to prove that the cycle class map
$c^4_{S\times \sh}$  is injective on $DCH^4(S\times
\sh)+\Q\,([o]\bt \iota ([o]))+\Q\,j_*q^*[o]$. Since we know by 
(\ref{sh}) and (\ref{gen}\ti a)) that it is injective on
$DCH^4(S\times \sh)$, this amounts to say that there is no
non-trivial relation
$$a\,[o]\bt \iota ([o]) + b\,j_*q^*[o]\in DH^8(S\times \sh)\ ,$$
with $a,b\in\Q$. \ind Suppose that such a relation holds.
 Let $\omega$ be a non-zero class in
$H^{2,0}(S)$; for any class
$\xi  
$ in $H^8(S\times \sh,\C)$  put $h(\xi 
):=(\pr_2)_*(\pr_1^*\omega
\cdot\xi )$. Since the product of $\omega $ with any algebraic
class in $H^2(S)$ is zero, $h$ is zero on $DH^8(S\times \sh)$.
Clearly  $h([o]\bt \iota ([o]))=0$, while $h(j_*q^*[o])=\pi
_*(p^*\omega \cdot q^*[o])=\iota (\omega )\iota ([o])$. This class
is non\-zero, for instance because  $\br{\iota (\omega )\iota
([o]),\iota (\bar
\omega )}=\br{\omega ,\bar\omega }>0$. 
\ind Thus $b=0$, and  our relation reduces to $[o]\bt \iota
([o])\in DH^8(S\times \sh)$. Since $ DH^8(S\times
\sh)=\sdir_{i+j=8}^{} DH^i(S)\bt DH ^j(\sh)$ (see Proposition
\ref{gen}\ti a)), this is equivalent to $\iota ([o])\in DH^4(\sh)$.
Thus the proof will be finished with the following lemma:
\th Lemma
\enonce The class $\iota([o])$ does
not belong to  $DH^4(S^{[2]})$.
\endth\label{io}\pv    We have
 $$\eqalign{  H^4(\ss)&=\varepsilon
^*H^4(S^2)\oplus i_*\eta^*H^2(S )\cr &=\Q\, p^*[o]\oplus
\Q\, q^*[o]\oplus (p^*H^2(S)\otimes q^*H^2(S))\oplus
i_*\eta^*H^2(S )\ .}$$ 
Taking the invariants under the involution 
of $\ss$ which exchanges the factors, we find
$$H^4(\sh)=\Q\, \iota ([o])\oplus \sym^2H^2(S)\oplus
 i'_*\eta^*H^2(S )\ ,$$
where  $\sym^2H^2(S)$ is identified with a
subspace of $H^4(\sh)$ by $\alpha \cdot \beta \mapsto \pi
_*(p^*\alpha
\cdot q^*\beta) $, and $i':=\pi \rond i$ is the natural embedding of
$E$ in $\sh$. Since
$CH^1(\sh)=$ \break$\iota (CH^1(S))\oplus\Q\cdot [E]$, the
subspace
$DH^4(\sh)$ is spanned by the classes 
$$\iota (\alpha)\iota (\beta)\ ,\quad
\iota (\alpha )\cdot [E]=2 i'_*\eta ^*\alpha \ ,\quad
[E]^2=-2\pi _*\varepsilon ^*[\Delta ] \quad \hbox{for }\ \alpha
,\beta\in CH^1(S)\ .$$
\ind Suppose that we have a relation
$$\iota ([o])=\sum_{i<j}m_{ij}\,\iota (\alpha _i)\iota (\alpha
_j)+i'_*\eta ^*\gamma +m\,\pi _*\varepsilon ^*[\Delta ]\quad
\hbox{in }\ H^4(\sh)\ ,$$
where $(\alpha _i)$ is a basis of $CH^1(S)$.
This gives in $H^4(\ss)$:
$$p^*[o]+q^*[o]=\sum_{i<j}m_{ij}\,(p^*\alpha
_i+q^*\alpha _i)(p^*\alpha _j+q^*\alpha _j)+2i_*\eta ^*\gamma
+2m\,\varepsilon ^*[\Delta ]\ .$$
Projecting onto the direct summand $i_*\eta^*H^2(S )$ of
$H^4(\ss)$ we find $i_*\eta
^*\gamma=0$. Multiplying by $p^*\omega $ and pushing forward
by $q$ we find as in the proof of (\ref{dh}) that all terms but
$\varepsilon ^*[\Delta ]$ give 0, so $m=0$. Finally 
the equality
$$p^*[o]+q^*[o]=\sum_{i<j}m_{ij}\,(p^*\alpha
_i+q^*\alpha _i)(p^*\alpha _j+q^*\alpha _j)$$ projected onto
$\sym^2H^2(S)$ gives $m_{ij}=0$ for all $i,j$. This achieves the
proof of the lemma and therefore of the Proposition.\cqfd 
\rem{Comments} A variation of this method can be used to prove
that the {\it generalized Kummer variety} $K_2$ associated to an
abelian surface $A$ [B1] has the weak splitting property; one must
replace $CH^1(A)$ by the subspace of {\it symmetric} divisor
classes. We leave the details to the reader. 
\ind We should point out, however, that even among symplectic
fourfolds these examples are quite particular. Indeed for each 
integer $g\ge 2$, the projective K3 surfaces of genus $g$ (that is,
embedded in $\P^g$ with degree $2g-2$) form an irreducible
19-dimensional family;   the corresponding family of Hilbert
schemes $S^{[2]}$  is contained in a 20-dimensional
irreducible family of projective symplectic manifolds (see [B1]). 
Since the
weak splitting property is not invariant under deformation, we do
not know whether it holds for the general member of such a
family. It would be interesting, in
particular, to check whether the property holds for the variety
of lines contained in a smooth cubic hypersurface in $\P^5$.

\vskip2cm
\centerline{ REFERENCES} \vglue15pt\baselineskip12.8pt
\def\num#1{\smallskip\item{\hbox to\parindent{\enskip [#1]\hfill}}}
\parindent=1.3cm 
\num{B1} A. {\pc BEAUVILLE}: {\sl 	Vari\'et\'es k\"ahl\'eriennes
dont la premi\`ere classe de Chern est nulle.}  J. of Diff. Geo\-metry
{\bf 18}, 755--782  (1983). 
\num{B2} A. {\pc BEAUVILLE}: {\sl 	Sur l'anneau de Chow d'une
vari\'et\'e ab\'elienne.} Math. Annalen {\bf 273}, 647--651 (1986).
\num{B-V} A. {\pc BEAUVILLE}, C. {\pc VOISIN}: {\sl  On the Chow
ring of a K3 surface}. J. Algebraic Geom., to appear.
\num{Bo} F. {\pc BOGOMOLOV}: {\sl On the cohomology ring of a
simple hyper-K\"ahler manifold (on the results of Verbitsky)}. Geom.
Funct. Anal. {\bf 6} (1996),  612--618.
\num{F} W. {\pc FULTON}: {\sl Intersection theory}. Ergebnisse
der Mathematik und ihrer Grenzgebiete (3)  {\bf 2}.
Springer-Verlag, Berlin (1984).
\num{J} U. {\pc JANNSEN}: {\sl Motivic sheaves and filtrations on
Chow groups}. Motives (Seattle, 1991), 245--302, Proc. Sympos.
Pure Math. {\bf 55} (Part 1), Amer. Math. Soc., Providence, RI, 1994.
\num{Jo} J.-P. {\pc JOUANOLOU}: {\sl Cohomologie de quelques
sch\'emas classiques et th\'eorie cohomologique des classes de
Chern}. SGA 5 (Cohomologie $l$\tx adique  et fonctions $L$),
Lecture Notes in Math. {\bf 589}, 282--350. Springer-Verlag, Berlin
(1977).
\num{H} D. {\pc HUYBRECHTS}: {\sl Compact hyper-K\"ahler
manifolds: basic results}. Invent. Math. {\bf 135} (1999),  63--113.
\num{L} M. {\pc LEHN}:  {\sl Chern classes of tautological sheaves  on
Hilbert schemes of points on surfaces}.  Invent. Math. {\bf 136} (1999), 
157--207.
\num{M} S. {\pc MUKAI}: {\sl Symplectic structure of the moduli
space of sheaves on an abelian or $K3$ surface}. Invent. Math.
{\bf 77} (1984), no. 1, 101--116.
\num{M-M} S. {\pc MORI}, S. {\pc MUKAI}: {\sl Classification of Fano
$3$-folds with} $B\sb{2}\geq 2$. Manu\-scri\-pta Math. {\bf 36}
(1981/82), 147--162. 
\num{S} J. {\pc SAWON}: {\sl Abelian fibred holomorphic
symplectic manifolds}. Turk. J. Math. {\bf 27}
(2003), 197-230. 
\num{W} J. {\pc WIERZBA}: {\sl Birational Geometry of Symplectic
4-folds}. Preprint (2002).

\vskip1cm
\font\eightrm=cmr8
\font\sixrm=cmr6
\def\pc#1{\eightrm#1\sixrm}
\hfill\vtop{\eightrm\hbox to 5cm{\hfill Arnaud {\pc BEAUVILLE}\hfill}
 \hbox to 5cm{\hfill Institut Universitaire de France\hfill}\vskip-2pt
\hbox to 5cm{\hfill \&\hfill}\vskip-2pt
 \hbox to 5cm{\hfill Laboratoire J.-A. Dieudonn\'e\hfill}
 \hbox to 5cm{\sixrm\hfill UMR 6621 du CNRS\hfill}
\hbox to 5cm{\hfill {\pc UNIVERSIT\'E DE}  {\pc NICE}\hfill}
\hbox to 5cm{\hfill  Parc Valrose\hfill}
\hbox to 5cm{\hfill F-06108 {\pc NICE} Cedex 02\hfill}}

\end